\documentclass[a4paper,12pt]{article}
\usepackage{amsmath}
\usepackage{amssymb}

\newtheorem{theorem}{Theorem}[section]

\newtheorem{corollary}[theorem]{Corollary}

\newtheorem{lemma}[theorem]{Lemma}

\newtheorem{proposition}[theorem]{Proposition}
\newtheorem{remark}[theorem]{Remark}

\newenvironment{proof}[1][Proof]{\noindent\textbf{#1.} }{\ \rule{0.5em}{0.5em}}
\input{tcilatex}
\begin{document}

\title{Some remarks on derivations in algebras of measurable operators}
\author{A.F. Ber, B. de Pagter and F.A. Sukochev\thanks{%
Research is partially supported by the Australian Research Council}}
\date{}
\maketitle

\begin{abstract}
This paper is concerned with derivations in algebras of (unbounded)
operators affiliated with a von Neumann algebra $\mathcal{M}$. Let $\mathcal{%
A}$ be one of the algebras of measurable operators, locally
measurable operators or, $\tau $-measurable operators. We
present a complete description of von Neumann algebras
$\mathcal{M}$ of type $I$ in terms of their central projections
such that every derivation in $\mathcal{A}$ is inner. It is
also shown that every derivation in the algebra $LS\left(
\mathcal{M}\right) $ of all locally measurable operators with
respect to a properly infinite von Neumann algebra
$\mathcal{M}$ vanishes on the center of $LS\left(
\mathcal{M}\right) $.
\end{abstract}

\section{Introduction}

Let $\mathcal{A}$ be a complex unitary algebra and $\mathcal{B}$ be an $%
\mathcal{A}$-bimodule. A linear map $\delta :\mathcal{A}\rightarrow \mathcal{%
B}$ is said to be a derivation if $\delta \left( xy\right) =\delta \left(
x\right) y+x\delta \left( y\right) $, $x,y\in \mathcal{A}$. It is an
important question whether for a given pair $\left( \mathcal{A},\mathcal{B}%
\right) $ every derivation $\delta $ is inner, that is, whether
there exists an element $y\in \mathcal{B}$ such that $\delta
=\delta _{y}$, where $\delta _{y}\left( x\right) =\left[
y,x\right] $, $x\in \mathcal{A}$. Here, $\left[ y,x\right]
=yx-xy$, the commutator of $x$ and $y$.

In the present paper, we consider the situation where $\mathcal{A}$ is a
subalgebra of $LS\left( \mathcal{M}\right) $, the algebra of all locally
measurable operators with respect to a von Neumann algebra $\mathcal{M}$,
and $\mathcal{B}\subseteq LS\left( \mathcal{M}\right) $ is a bi-module over $%
\mathcal{A}$ (with respect to the multiplication in $LS\left( \mathcal{M}%
\right) $). In Section \ref{SectTypeI}, a characterization of type $I$ von
Neumann algebras $\mathcal{M}$ is presented for which all derivations $%
\delta :\mathcal{A}\rightarrow \mathcal{A}$ are inner, where $\mathcal{A}$
is any of the algebras $LS\left( \mathcal{M}\right) $, $S\left( \mathcal{M}%
\right) $ or $S\left( \tau \right) $; here $\tau $ denotes a semi-finite,
faithful normal trace on $\mathcal{M}$ (for definitions see Section \ref%
{SectPre} and the references therein). In Section \ref{SectZlin}, it is
shown that, for every properly infinite von Neumann algebra $\mathcal{M}$,
any derivation $\delta :\mathcal{A}\rightarrow LS\left( \mathcal{M}\right) $
vanishes on the center $\mathcal{Z}\left( \mathcal{A}\right) $, where $%
\mathcal{A}$ is any $\ast $-subalgebra of $LS\left( \mathcal{M}\right) $
such that $\mathcal{M}\subseteq \mathcal{A}$. These results extend and
complement the classical results in the theory of derivations in von Neumann
algebras (see \cite{Dix}, \cite{Sak}). It is also a classical result that
any derivation $\delta :\mathcal{M}\rightarrow \mathcal{M}$ is inner and
that there exists $a\in \mathcal{M}$ such that $\delta =\delta _{a}$ and $%
\left\Vert a\right\Vert \leq \left\Vert \delta \right\Vert $ (and, actually,
$a$ may be chosen such that $2\left\Vert a\right\Vert \leq \left\Vert \delta
\right\Vert $). The algebras $LS\left( \mathcal{M}\right) $, $S\left(
\mathcal{M}\right) $ and $S\left( \tau \right) $ are not normed and so, the
norm of a derivation in such algebra does not make sense. However, in
Section \ref{SectEst}, we provide interesting operator inequalities which
may be considered as substitutes for norm inequalities, and provide
additional insight even in the classical normed case.

\section{Some preliminaries\label{SectPre}}

Let $\mathcal{A}$ be a unital algebra. The center of $\mathcal{A}$ is
denoted by $\mathcal{Z}\left( \mathcal{A}\right) $, that is, $\mathcal{Z}%
\left( \mathcal{A}\right) =\left\{ a\in \mathcal{A}:ax=xa\ \forall x\in
\mathcal{A}\right\} $. Furthermore, $\mathcal{I}\left( \mathcal{A}\right) $
denotes the set of all idempotents in $\mathcal{A}$. Let $\mathcal{B}$ be a
bi-module over $\mathcal{A}$, where it is assumed that $ax=xa$ for all $a\in
\mathcal{Z}\left( \mathcal{A}\right) $ and $x\in \mathcal{B}$. A linear map $%
\delta :\mathcal{A}\rightarrow \mathcal{B}$ is called a \textit{derivation}
if $\delta \left( xy\right) =\delta \left( x\right) y+x\delta \left(
y\right) $, $x,y\in \mathcal{A}$. In particular, if $\mathcal{B}=\mathcal{A}$%
, then such a map is called a \textit{derivation in }$\mathcal{A}$. If $b\in
\mathcal{B}$, then the map $\delta _{b}:\mathcal{A}\rightarrow \mathcal{B}$,
given by $\delta _{b}\left( x\right) =bx-xb$, $x\in \mathcal{A}$, is a
derivation. Such a derivation is called \textit{inner}. Furthermore, a
derivation $\delta :\mathcal{A}\rightarrow \mathcal{B}$ is called $\mathcal{Z%
}$\textit{-linear} if $\delta \left( ax\right) =a\delta \left( x\right) $
for all $a\in \mathcal{Z}\left( \mathcal{A}\right) $ and $x\in \mathcal{A}$.
Evidently, $\delta $ is $\mathcal{Z}$-linear if and only if $\delta \left(
a\right) =0$ for all $a\in \mathcal{Z}\left( \mathcal{A}\right) $. We also
recall the following simple properties:

\begin{enumerate}
\item[($\protect\alpha $).] If $a\in \mathcal{Z}\left( \mathcal{A}\right) $,
then $x\delta \left( a\right) =\delta \left( a\right) x$ for all $x\in
\mathcal{A}$. Indeed,
\begin{equation*}
x\delta \left( a\right) =\delta \left( xa\right) -\delta \left( x\right)
a=\delta \left( ax\right) -a\delta \left( x\right) =\delta \left( a\right) x.
\end{equation*}%
In particular, if $\mathcal{B}=\mathcal{A}$, then $\delta \left( a\right)
\in \mathcal{Z}\left( \mathcal{A}\right) $ whenever $a\in \mathcal{Z}\left(
\mathcal{A}\right) $.

\item[($\protect\beta $).] If $p\in \mathcal{I}\left( \mathcal{Z}\left(
\mathcal{A}\right) \right) $, then $\delta \left( p\right) =0$ and $\delta
\left( px\right) =p\delta \left( x\right) $ for all $x\in \mathcal{A}$.

\item[($\protect\gamma $).] If $p\in \mathcal{I}\left( \mathcal{A}\right) $
and $x\in \mathcal{A}$ such that $xp=px$, then $p\delta \left( x\right)
p=p\delta \left( xp\right) p$. Indeed, it is easy to see that $p\delta
\left( p\right) p=0$. Therefore, multiplying the identity $\delta \left(
xp\right) =\delta \left( x\right) p+x\delta \left( p\right) $ left and right
by $p$, the result follows.
\end{enumerate}

In this paper, $\mathcal{M}$ denotes a von Neumann algebra on a Hilbert
space $H$. The $\ast $-algebras of all measurable operators and all locally
measurable operators, with respect to $\mathcal{M}$, are denoted by $S\left(
\mathcal{M}\right) $ and $LS\left( \mathcal{M}\right) $, respectively.
Furthermore, if $\tau $ is a faithful normal semi-finite trace on $\mathcal{M%
}$, then the $\ast $-algebra of all $\tau $-measurable operators is denoted
by $S(\mathcal{M},\tau )=S\left( \tau \right) $ (for the definitions, we
refer the reader to e.g. \cite{Se}, \cite{San}, \cite{Ye} and \cite{Ne}).

The set of all (self-adjoint) projections in $\mathcal{M}$ is denoted by $%
P\left( \mathcal{M}\right) $. A non-zero projection $p\in \mathcal{P}(%
\mathcal{M})$ is said to be an $\emph{atom}$ (or, \textit{minimal}) if $%
0<q\leq p$ in $\mathcal{P}(\mathcal{M})$ implies that $q=p$. A non-zero
projection $p\in \mathcal{P}(\mathcal{M})$ is called $\emph{discrete}$ if $p$
is a supremum of atoms in $P\left( \mathcal{M}\right) $. A non-zero
projection $p\in \mathcal{P}(\mathcal{M})$ is said to be $\emph{continuous}$
if there are no atoms $q\in \mathcal{P}(\mathcal{M})$ satisfying $0<q\leq p$%
. If a projection $0\neq p\in P\left( \mathcal{M}\right)$ is not
discrete, then it follows via a simple maximality argument that
$p$ dominates a continuous projection in $P\left(
\mathcal{M}\right) $.

\section{Von Neumann algebras of type $I$\label{SectTypeI}}

If $\left( X,\Sigma ,\mu \right) $ is a Maharam measure space (that is, $\mu
$ has the finite subset property and the measure algebra is a complete
Boolean algebra), then the algebra $\mathcal{M}=L_{\infty }\left( \mu
\right) $ is a von Neumann algebra (acting via multiplication on the Hilbert
space $H=L_{2}\left( \mu \right) $). In this case we have $LS\left( \mathcal{%
M}\right) =S\left( \mathcal{M}\right) =L_{0}\left( \mu \right) $, where $%
L_{0}\left( \mu \right) $ is the algebra of all (equivalence classes of)
complex valued $\mu $-measurable functions. Denoting by $\mu $ also the
corresponding trace on $\mathcal{M}$ (that is, $\mu \left( f\right)
=\int_{X}fd\mu $, $f\in \mathcal{M}^{+}$), it should be observed that, in
general, $S\left( \mu \right) $ is a proper subalgebra of $L_{0}\left( \mu
\right) $. The following observation will be used.

\begin{lemma}
\label{Lem01}If $\left( X,\Sigma ,\mu \right) $ is a Maharam measure space
and $M_{n}$ ($n\in \mathbb{N}$) denotes the von Neumann algebra of all
complex $n\times n$-matrices (equipped with the standard trace $\limfunc{tr}%
_{n}$), then $LS\left( L_{\infty }\left( \mu \right) \otimes M_{n}\right)
=S\left( L_{\infty }\left( \mu \right) \otimes M_{n}\right) \cong
L_{0}\left( \mu \right) \otimes M_{n}$ and $S\left( \mu \otimes \limfunc{tr}%
_{n}\right) \cong S\left( \mu \right) \otimes M_{n}$.
\end{lemma}

In the proof of Theorem \ref{t1} the following result will be used.

\begin{theorem}[\protect\cite{BCS2}, Theorem 3.4]
\label{Thm2}If $\left( X,\Sigma ,\mu \right) $ is a Maharam
measure space, then the following two statements are
equivalent:

\begin{enumerate}
\item[(i).] there exists a non-zero derivation in the algebra $L_{0}\left(
\mu \right) $;

\item[(ii).] the measure algebra of $\left( X,\Sigma ,\mu \right) $ is not
atomic.
\end{enumerate}
\end{theorem}

In the case when $\mathcal{M}$ is a von Neumann algebra of type
$I$ and $\mathcal{A}$ is either $LS\left( \mathcal{M}\right)$,
or $S\left( \mathcal{M}\right)$, or $S\left( \mathcal{M},\tau
\right)$, the description of all derivations in $\mathcal{A}$
is given in \cite{AAK}. The following theorem gives a necessary
and sufficient condition on the algebra $\mathcal{M}$
guaranteeing that any derivation on  $\mathcal{A}$ is inner.
The proof is based on the description of the algebra $LS\left(
\mathcal{M}\right)$ as an algebra of $B(H)$-valued functions
established in \cite{BdPS}.

\begin{theorem}
\label{t1} Let $\mathcal{M}$ be von Neumann algebra of type $I$ with a
separable pre-dual $\mathcal{M}_{\ast }$ and let $\tau $ be a faithful
normal semi-finite trace on $\mathcal{M}$. If $\mathcal{A}$ denotes one of
the algebras $LS(\mathcal{M})$, $S(\mathcal{M})$ or $S(\mathcal{M},\tau )$,
then the following statements are equivalent:

\begin{enumerate}
\item[(i).] every derivation on $\mathcal{A}$ is inner;

\item[(ii).] every central projection in $\mathcal{M}$ is either infinite or
discrete.
\end{enumerate}
\end{theorem}

\begin{proof}
Since $\mathcal{M}$ is of type $I$ and the space $\mathcal{M}_{\ast }$ is
separable, there exists a (unique) pairwise orthogonal family $%
\{z_{n}\}_{n=0}^{\infty }\subset P(\mathcal{Z}(\mathcal{M})))$ such that $%
\sum_{n=0}^{\infty }z_{n}=1$ and
\begin{eqnarray*}
\mathcal{M}z_{n} &\cong &L_{\infty }\left( \mu _{n}\right) \otimes M_{n},\ \
n>0, \\
\mathcal{M}z_{0} &\cong &L_{\infty }(\mu _{0})\overline{\otimes }B(H_{0}),
\end{eqnarray*}%
where $M_{n}=B\left( H_{n}\right) $ is the algebra of $n\times n$ complex
matrices (that is, $H_{n}=\mathbb{C}^{n}$), $H_{0}=L_{2}(0,1)$ and $\mu _{n}$
are $\sigma $-finite and separable measures on $\left( X_{n},\Sigma
_{n}\right) $ for all $n\geq 0$. If $z_{n}=0$ for some $n\geq 0$, then the
corresponding space $X_{n}=\emptyset $ (see e.g. \cite{KR}, Section 6.6).

(i)$\Rightarrow $(ii). Let $0\neq p\in P(\mathcal{Z}(\mathcal{M}))$ be given
and suppose that $p$ is finite and not discrete. We shall derive a
contradiction. As observed above, $p$ dominates a continuous projection in $%
P\left( \mathcal{Z}\left( \mathcal{M}\right) \right) $ and so, we may assume
that $p$ is finite and continuous. If $pz_{0}\neq 0$, then $pz_{0}$ is a
central projection in $L_{\infty }(\mu _{0})\overline{\otimes }B(H_{0})$ and
hence, $pz_{0}=\chi _{A}\otimes \mathbf{1}_{H_{0}}$ for some $A\in \Sigma
_{0}$ with $\mu _{0}\left( A\right) >0$. Since $H_{0}$ is infinite
dimensional, this implies that $pz_{0}$ is infinite. Since $p$ is finite,
this is impossible. Consequently, $pz_{0}=0$. This implies that there exists
$n>0$ such that $pz_{n}\neq 0$. Evidently, $\mathcal{M}pz_{n}\cong L_{\infty
}(\mu )\otimes M_{n}$, where the measure $\mu $ on $\left( X,\Sigma \right) $
is atomless and $\sigma $-finite. The restriction of $\tau $ to $\mathcal{M}%
pz_{n}$ is a faithful normal semifinite trace and so, it follows from \cite%
{Se} that this restriction is given by $\nu \otimes \func{tr}_{n}$, where $%
\func{tr}_{n}$ is the standard trace on $M_{n}$ and $\nu $ is a $\sigma $%
-finite measure on $\left( X,\Sigma \right) $ which is equivalent with $\mu $
(so, $\nu \left( f\right) =\int_{X}fd\nu $ for all $f\in L_{\infty }\left(
\mu \right) ^{+}=L_{\infty }\left( \nu \right) ^{+}$. If $A\in \Sigma $ is
such that $0<\nu \left( A\right) <\infty $, then the projection $p_{1}=\chi
_{A}\otimes \mathbf{1}_{M_{n}}$ satisfies $p_{1}\in \mathcal{Z}(\mathcal{M}%
),\ 0<p_{1}\leq pz_{n}$, and $\tau (p_{1})<\infty $. Consequently, $\mathcal{%
A}p_{1}\cong L_{0}\left( \nu _{1}\right) \otimes M_{n}$, where $\nu _{1}$ is
the restriction to the set $A$.

It follows from Theorem \ref{Thm2}, that there exists a derivation $\delta
_{0}\neq 0$ on $L_{0}\left( \nu _{1}\right) $. Define the derivation $\delta
_{1}$ on $L_{0}\left( \nu _{1}\right) \otimes M_{n}$ by setting
\begin{equation*}
\delta _{1}(a\otimes b)=\delta _{0}(a)\otimes b,\ \ \ a\in L_{0}\left( \nu
_{1}\right) ,b\in M_{n}.
\end{equation*}%
We claim that is not inner. Indeed, if $a\in L_{0}\left( \nu _{1}\right) $
and $\delta _{0}(a)\neq 0$, then $\delta _{1}(a\otimes \mathbf{1}%
_{M_{n}})=\delta _{0}(a)\otimes \mathbf{1}_{M_{n}}\neq 0$. The element $%
a\otimes \mathbf{1}_{M_{n}}$ belongs to the center of the algebra $%
L_{0}\left( \nu _{1}\right) \otimes M_{n}$ and any inner derivation vanishes
on the center. This proves our claim.

To the derivation $\delta _{1}$ on $L_{0}\left( \nu _{1}\right) \otimes
M_{n} $ corresponds a derivation $\delta _{2}$ on $\mathcal{A}pz_{n}$, which
is obviously not inner. Therefore, the derivation $\delta $ on $\mathcal{A}$%
, defined by setting $\delta (x):=\delta _{2}(xpz_{n})$, $x\in \mathcal{A}$,
is also not inner. This contradicts assumption (i) and so, we may conclude
that (i) implies (ii).

(ii)$\Rightarrow $(i). Let $\delta $ be a derivation on the algebra $%
\mathcal{A}$. By ($\beta $) in Section \ref{SectPre}, it is clear that $%
\delta $ leaves $\mathcal{A}z_{n}$ invariant for each $n\geq 0$. Let $\delta
_{n}$ be the restriction of $\delta $ to $\mathcal{A}z_{n}$. If $n>0$, then
the von Neumann algebra $\mathcal{M}z_{n}$ is finite and so, the projection $%
z_{n}$ is also finite. Therefore, by hypothesis, $z_{n}$ is discrete for $%
n>0 $. Hence, the measure spaces $(X_{n},\Sigma _{n},\mu _{n})$ are atomic
for all $n>0$.

We first consider the cases $\mathcal{A}=LS(\mathcal{M})$ or $S(\mathcal{M}%
). $ If $n>0$, then $\mathcal{A}z_{n}\cong L_{0}(\mu _{n})\otimes M_{n}$
(see Lemma \ref{Lem01}) and
\begin{equation*}
\mathcal{Z}(L_{0}(\mu _{n})\otimes M_{n})=L_{0}(\mu _{n})\otimes \mathbb{C}%
\mathbf{1}\cong L_{0}(\mu _{n}),\ \ \ n>0.
\end{equation*}%
Since $\delta _{n}$ maps $\mathcal{Z}\left( \mathcal{A}z_{n}\right) $ into
itself, it follows from from Theorem \ref{Thm2} $\delta _{n}=0$ on $\mathcal{%
Z}\left( \mathcal{A}z_{n}\right) $, that is, $\delta _{n}$ is $\mathcal{Z}$%
-linear. Therefore, by \cite{BdPS}, Corollary 6.8, each $\delta _{n}$, $n>0$
is inner. Furthermore, since
\begin{equation*}
LS(\mathcal{M})z_{0}\cong LS(L_{\infty }(\mu _{0})\overline{\otimes }%
B(H_{0})),\ \ \ \ S(\mathcal{M})z_{0}\cong S(L_{\infty }(\mu _{0})\overline{%
\otimes }B(H_{0})),
\end{equation*}%
where $H_{0}$ is infinite dimensional and separable, we infer from \cite%
{BdPS}, Corollary 6.19, that $\delta _{0}$ is inner.

Now consider the case $\mathcal{A}=S(\tau )$. The restriction $\tau _{n}$ of
the trace $\tau $ to $\mathcal{M}z_{n}=L_{\infty }\left( \mu _{n}\right)
\overline{\otimes }B\left( H_{n}\right) $ is a faithful normal semifinite
trace for every $n\geq 0$. It follows from \cite{Se} that $\tau _{n}=\nu
_{n}\otimes \limfunc{tr}_{n}$, where $\nu _{n}$ is a measure on $\left(
X,\Sigma _{n}\right) $ which is equivalent with $\mu _{n}$ and $\limfunc{tr}%
_{n}$ denotes the standard trace on $B\left( H_{n}\right) $. Consequently, $%
\mathcal{A}z_{n}=S\left( L_{\infty }\left( \mu _{n}\right) \overline{\otimes
}B\left( H_{n}\right) ,\nu _{n}\otimes \limfunc{tr}_{n}\right) $ for all $%
n\geq 0$. If $n>0$, then $B\left( H_{n}\right) =M_{n}$ and hence, by Lemma %
\ref{Lem01},
\begin{equation*}
\mathcal{A}z_{n}=S\left( \nu _{n}\right) \otimes M_{n},\ \ \ \mathcal{Z}%
\left( \mathcal{A}z_{n}\right) \cong S\left( \nu _{n}\right) .
\end{equation*}%
Since the measure $\nu _{n}$ is discrete and $\delta _{n}$ maps $\mathcal{Z}%
\left( \mathcal{A}z_{n}\right) $ into itself, it follows that $\delta _{n}$
vanishes on $\mathcal{Z}\left( \mathcal{A}z_{n}\right) $. Hence, $\delta
_{n} $ is $\mathcal{Z}$-linear and so, it follows from \cite{BdPS},
Corollary 6.18, that $\delta _{n}$ is inner for every $n>0$. For $n=0$ it is
an immediate consequence of \cite{BdPS}, Corollary 6.19, that $\delta _{0}$
is inner.

We thus, have shown that in each of the three cases there exist $a_{n}\in
\mathcal{A}z_{n}$ such that $\delta _{n}=\left[ a_{n},.\right] $ for all $%
n\geq 0$. Define the operator $a\in LS\left( \mathcal{M}\right) $ by setting
$a=\sum_{n=0}^{\infty }a_{n}\in LS(\mathcal{M})$. It is easily verified that
$\delta \left( x\right) =\left[ a,x\right] $, $x\in \mathcal{A}$. Since $S(%
\mathcal{M})$ and $S(\tau )$ are absolutely solid $\ast $-subalgebras in $LS(%
\mathcal{M})$, containing the algebra $\mathcal{M}$, we may apply \cite{BdPS}%
, Proposition 6.17, which guarantees that $\delta $ is an inner derivation
on $\mathcal{A}$ in all three cases.
\end{proof}

\section{$\mathcal{Z}$-linearity\label{SectZlin}}

In this section, $\mathcal{M}$ denotes a von Neumann algebra on a Hilbert
space $H$. For the proof of the main result (Proposition \ref{AProp01}), we
need some preparations. Using the center-valued trace on the reduced finite
von\ Neumann algebra $e\mathcal{M}e$, where $e=e_{1}\vee p$, the proof of
the first of the following lemmas is straightforward.

\begin{lemma}
\label{ALem02} (see e.g. \cite{BdPS}, Lemma 6.11). If $\left\{ e_{n}\right\}
_{n=1}^{\infty }$ is a sequence of finite projections in $\mathcal{M}$ such
that $e_{n}\downarrow 0$ and if $p\in P\left( \mathcal{M}\right) $ is such
that $p\precsim e_{n}$ for all $n$, then $p=0$.
\end{lemma}

\begin{lemma}
\label{ALem01}If $a\in LS\left( \mathcal{M}\right) ^{+}$ and $q\in P\left(
\mathcal{M}\right) $ are such that $qaq\geq \lambda q$ for some $0<\lambda
\in \mathbb{R}$, then $q\precsim e^{a}\left( \mu ,\infty \right) $ for all $%
0\leq \mu <\lambda $.
\end{lemma}

\begin{proof}
Given $0\leq \mu <\lambda $, define $p\in P\left( \mathcal{M}\right) $ by $%
p=q\wedge e^{a}\left[ 0,\mu \right] $. It follows that
\begin{eqnarray*}
\mu p &=&p\left( \mu e^{a}\left[ 0,\mu \right] \right) p\geq pae^{a}\left[
0,\mu \right] p=pap=p\left( qaq\right) p \\
&\geq &p\left( \lambda q\right) p=\lambda p.
\end{eqnarray*}%
Since $\mu <\lambda $, this implies that $p=0$. Hence, it follows from
Kaplansky's formula that
\begin{equation*}
q=q-q\wedge e^{a}\left[ 0,\mu \right] \sim q\vee e^{a}\left[ 0,\mu \right]
-e^{a}\left[ 0,\mu \right] \leq e^{a}\left( \mu ,\infty \right) ,
\end{equation*}%
that is, $q\precsim e^{a}\left( \mu ,\infty \right) $.
\end{proof}

\begin{lemma}
\label{ALem06}Suppose that $a\in LS\left( \mathcal{M}\right) ^{+}$. If $%
\left\{ q_{n}\right\} _{n=1}^{\infty }$ is a sequence in $P\left( \mathcal{M}%
\right) $ and $q_{0}\in P\left( \mathcal{M}\right) $ such that $q_{n}\sim
q_{0}$ and $q_{n}aq_{n}\geq nq_{n}$ for all $n\in \mathbb{N}$, then $q_{0}=0$%
.
\end{lemma}

\begin{proof}
It will be assumed first that $a\in S\left( \mathcal{M}\right) ^{+}$. It
follows from Lemma \ref{ALem01} that $q_{n}\precsim e^{a}\left( n-1,\infty
\right) $ and so, $q_{0}\precsim e^{a}\left( n-1,\infty \right) $ for all $%
n\in \mathbb{N}$. Since $a\in S\left( \mathcal{M}\right) ^{+}$, there exists
$n_{0}\in \mathbb{N}$ such that $e^{a}\left( n-1,\infty \right) $ is finite
for all $n\geq n_{0}$. Evidently, $e^{a}\left( n-1,\infty \right) \downarrow
0$ and hence, Lemma \ref{ALem02} implies that $q_{0}=0$.

Now, assume that $a\in LS\left( \mathcal{M}\right) ^{+}$. By definition,
there exists a sequence $\left\{ z_{k}\right\} _{k=1}^{\infty }$ in $P\left(
\mathcal{Z}\left( \mathcal{M}\right) \right) $ such that $z_{k}\uparrow
\mathbf{1}$ and $az_{k}\in S\left( \mathcal{M}\right) $ for all $k$. It
follows from $q_{n}aq_{n}\geq nq_{n}$ that $q_{n}\left( az_{k}\right)
q_{n}\geq nq_{n}z_{k}$ for all $n\in \mathbb{N}$. Since $q_{n}z_{k}\sim
q_{0}z_{k}$ for all $n$, the first part of the proof implies that $%
q_{0}z_{k}=0$ for all $k$. Since $q_{0}z_{k}=q_{0}z_{k}q_{0}\uparrow
_{k}q_{0}$, it is now clear that $q_{0}=0$. \bigskip
\end{proof}

In the next theorem, we assume that $\mathcal{A}$ is a $\ast
$-subalgebra of $LS\left( \mathcal{M}\right) $ such that
$\mathcal{M}\subseteq \mathcal{A}$. We consider $LS\left(
\mathcal{M}\right) $ as a bi-module over $\mathcal{A}$ (with
respect to the multiplication in $LS\left( \mathcal{M}\right)
$). It
should be noted that $ax=xa$ for all $a\in \mathcal{Z}\left( \mathcal{A}%
\right) $ and $x\in LS\left( \mathcal{M}\right) $. Indeed, if $a\in \mathcal{%
Z}\left( \mathcal{A}\right) $, then $ax=xa$ for all $x\in \mathcal{M}$ and
hence, the same holds for all $x\in LS\left( \mathcal{M}\right) $. Now, we
are in a position to prove the following result.

\begin{theorem}
\label{AProp01} Suppose that $\mathcal{M}$ is properly infinite,
that is, every non-zero projection in $P\left( \mathcal{Z}\left(
\mathcal{M}\right)
\right) $ is infinite (with respect to $\mathcal{M}$). Let $\mathcal{A}%
\subseteq LS\left( \mathcal{M}\right) $ be a $\ast $-subalgebra such that $%
\mathcal{M}\subseteq \mathcal{A}$. If $\delta :\mathcal{A}\rightarrow
LS\left( \mathcal{M}\right) $ is a derivation, then $\delta $ is $\mathcal{Z}
$-linear (that is, $\delta \left( ax\right) =a\delta \left( x\right) $ for
all $x\in \mathcal{A}$ and $a\in \mathcal{Z}\left( \mathcal{A}\right) $).
\end{theorem}

\begin{proof}
We start the proof by observing that $\delta \left( a\right) \in \mathcal{Z}%
\left( LS\left( \mathcal{M}\right) \right) $ whenever $a\in \mathcal{Z}%
\left( \mathcal{A}\right) $. Indeed, it follows from ($\alpha $) in Section %
\ref{SectPre} that $x\delta \left( a\right) =\delta \left( a\right) x$ for
all $x\in \mathcal{A}$ and so, in particular for all $x\in \mathcal{M}$.
From this it follows easily that $\delta \left( a\right) \in \mathcal{Z}%
\left( LS\left( \mathcal{M}\right) \right) $.

To show that $\delta $ is $\mathcal{Z}$-linear, is clearly sufficient to
prove that $\delta \left( a\right) =0$ for all $a\in \mathcal{Z}\left(
\mathcal{A}\right) $. Furthermore, it may be assumed, without loss of
generality, that $\delta $ is self-adjoint (that is, $\delta \left( x^{\ast
}\right) =\delta \left( x\right) ^{\ast }$ for all $x\in \mathcal{A}$).

Suppose that $\delta \neq 0$ on $\mathcal{Z}\left( \mathcal{A}\right) $.
This implies that there exists $a=a^{\ast }\in \mathcal{Z}\left( \mathcal{A}%
\right) $ such that $\delta \left( a\right) \neq 0$. Without loss of
generality, it may be assumed that $\delta \left( a\right) ^{+}\neq 0$. Let $%
s$ be the support projection of $\delta \left( a\right) ^{+}$. Since $\delta
\left( a\right) \in \mathcal{Z}\left( LS\left( \mathcal{M}\right) \right) $,
it follows that $s\in P\left( \mathcal{Z}\left( \mathcal{M}\right) \right) $
and so, by ($\beta $) in Section \ref{SectPre}, $\delta \left( sa\right)
=s\delta \left( a\right) =\delta \left( a\right) ^{+}>0$. Consequently,
replacing $a$ by $sa$, it may be assumed, without loss of generality, that $%
a=a^{\ast }\in \mathcal{Z}\left( \mathcal{A}\right) $ satisfies $\delta
\left( a\right) >0$. Multiplying $a$ by a positive scalar, it may
furthermore be assumed that $\delta \left( a\right) \geq p_{0}$ for some $%
0\neq p_{0}\in P\left( \mathcal{Z}\left( \mathcal{M}\right) \right) $.

For $n\in \mathbb{N}$, let $f_{n}:\mathbb{R}\rightarrow \mathbb{R}$ be
defined by%
\begin{equation*}
f_{n}\left( \lambda \right) =\sum_{k=-\infty }^{\infty }\left( n\lambda
-k\right) \chi _{\left( \frac{k}{n},\frac{k+1}{n}\right] },\ \ \ \lambda \in
\mathbb{R}.
\end{equation*}%
Since $0\leq f_{n}\left( \lambda \right) \leq \mathbf{1}$, $\lambda \in
\mathbb{R}$, it is clear that $0\leq f_{n}\left( a\right) \leq \mathbf{1}$
in $\mathcal{Z}\left( \mathcal{M}\right) $. Since $\mathcal{M}\subseteq
\mathcal{A}$, this implies, in particular, that $f_{n}\left( a\right) \in
\mathcal{Z}\left( \mathcal{A}\right) $. Furthermore,
\begin{equation*}
0\leq \left( n\lambda -f_{n}\left( \lambda \right) \right) \chi _{\left(
\frac{k}{n},\frac{k+1}{n}\right] }=k\chi _{\left( \frac{k}{n},\frac{k+1}{n}%
\right] }
\end{equation*}%
and so, $\left( na-f_{n}\left( a\right) \right) e^{a}\left( \frac{k}{n},%
\frac{k+1}{n}\right] =ke^{a}\left( \frac{k}{n},\frac{k+1}{n}\right] $.
Hence, by ($\beta $) in Section \ref{SectPre},
\begin{equation*}
\left( n\delta \left( a\right) -\delta \left( f_{n}\left( a\right) \right)
\right) e^{a}\left( \frac{k}{n},\frac{k+1}{n}\right] =0,\ \ \ k\in \mathbb{Z}%
.
\end{equation*}%
Consequently, $n\delta \left( a\right) =\delta \left( f_{n}\left( a\right)
\right) $ for all $n\in \mathbb{N}$. Defining $a_{n}=f_{n}\left( a\right) $,
it follows that $0\leq a_{n}\leq \mathbf{1}$ and $\delta \left( a_{n}\right)
=n\delta \left( a\right) \geq np_{0}$ for all $n\in \mathbb{N}$.

Since $0\neq p_{0}\in P\left( \mathcal{Z}\left( \mathcal{M}\right) \right) $
is an infinite projection in $\mathcal{M}$, there exists a sequence $\left\{
p_{n}\right\} _{n=1}^{\infty }$ of mutually orthogonal projections in $%
P\left( \mathcal{M}\right) $ such that $p_{n}\leq p_{0}$ and
$p_{n}\sim p_{0} $ for all $n\in \mathbb{N}$ (\cite{KR}, Lemma
6.3.3).
Define $b\in \mathcal{M}$ by setting $%
b=\sum_{n=1}^{\infty }a_{n}p_{n}$ (as a series converging in the
strong operator topology). It should be observed that
$p_{n}b=bp_{n}=a_{n}p_{n}$ for all $n\in \mathbb{N}$. Using
($\gamma $) in Section \ref{SectPre}, twice, it follows that
\begin{equation*}
p_{n}\delta \left( b\right) p_{n}=p_{n}\delta \left( bp_{n}\right)
p_{n}=p_{n}\delta \left( a_{n}p_{n}\right) p_{n}=p_{n}\delta
\left( a_{n}\right) p_{n}
\end{equation*}%
and so,
\begin{equation*}
p_{n}\left\vert \delta \left( b\right) \right\vert p_{n}\geq
p_{n}\delta \left( b\right) p_{n}=p_{n}\delta \left( a_{n}\right)
p_{n}\geq p_{n}\left( np_{0}\right) p_{n}=np_{n}
\end{equation*}%
for all $n\in \mathbb{N}$. Lemma \ref{ALem06} implies that $p_{0}=0$. This
is a contradiction, by which the proof is complete.\bigskip
\end{proof}

The following result is an immediate consequence of the above
theorem.

\begin{corollary}
Suppose that $\mathcal{M}$ is properly infinite. Let $\mathcal{A}\subseteq
LS\left( \mathcal{M}\right) $ be a $\ast $-subalgebra such that $\mathcal{M}%
\subseteq \mathcal{A}$ and $\mathcal{B}\subseteq LS\left( \mathcal{M}\right)
$ be an $\mathcal{A}$-bimodule. If $\delta :\mathcal{A}\rightarrow \mathcal{B%
}$ is a derivation, then $\delta $ is $\mathcal{Z}$-linear.
\end{corollary}

\section{Operator estimates\label{SectEst}}

If $\delta $ is a derivation on an algebra $\mathcal{A}$ which is inner,
then, by definition, there exists $a\in \mathcal{A}$ such that $\delta
\left( a\right) =\left[ a,x\right] $, $x\in \mathcal{A}$. The element $a$ is
not uniquely determined by $\delta $. Indeed, if $b\in \mathcal{Z}\left(
\mathcal{A}\right) $, then $a-b$ induces the same derivation $\delta $. The
question arises whether there exist any "good" elements $a$ in $\mathcal{A}$
which induce $\delta $. By way of example, if $\mathcal{M}$ is a von Neumann
algebra, then every derivation $\delta $ on $\mathcal{M}$ is inner and there
exists $a\in \mathcal{M}$ such that $\delta \left( x\right) =\left[ a,x%
\right] $, $x\in \mathcal{A}$, and $\left\Vert a\right\Vert \leq
\left\Vert \delta \right\Vert $ (see e.g. \cite{Dix}, Chapter 9).
Moreover, if $\delta $ is self-adjoint, then $a$ may even be
chosen such that $2\left\Vert a\right\Vert \leq \left\Vert \delta
\right\Vert $ (see \cite{BR}, Corollary 3.2.47). If $\mathcal{A%
}=LS\left( \mathcal{M}\right) $, then the algebra $\mathcal{A}$ is not
normed (in general) and so, norm estimates do not make sense. However, one
may hope to replace norm estimates by operator estimates. By way of example,
in \cite{BdPS}, Lemma 6.16, the following result has been obtained.

\begin{proposition}
\label{Prop02}If $a\in LS_{h}\left( \mathcal{M}\right) $, then there exists
a partial isometry $v\in \mathcal{M}$ and an operator $b\in \mathcal{Z}%
\left( LS_{h}\left( \mathcal{M}\right) \right) $ such that $\left\vert \left[
a,v\right] \right\vert +\mathbf{1}\geq \left\vert a-b\right\vert $.
\end{proposition}

As an immediate application it follows that if $\mathcal{A}$ is an
absolutely solid $\ast $-subalgebra of $LS\left( \mathcal{M}\right) $ such
that $\mathcal{M}\subseteq \mathcal{A}$, and if $\delta $ is a derivation in
$\mathcal{A}$ which is given by $\delta \left( x\right) =\left[ a,x\right] $%
, $x\in \mathcal{A}$, for some $a\in \mathcal{A}$, then $a$ may be chosen
from $\mathcal{A}$.

The purpose of the present section is to obtain some sharpening of
Proposition \ref{Prop02}. For this purpose, it will be convenient
to recall some of the ingredients in the proof of this
proposition, in particular, the
construction of the element $b\in \mathcal{Z}\left( LS_{h}\left( \mathcal{M}%
\right) \right) $. If $e,f\in P\left( \mathcal{M}\right) $, then the central
projection $z\left( e,f\right) $ in $\mathcal{M}$ is defined by setting
\begin{equation}
z\left( e,f\right) =\bigvee \left\{ p\in P\left( \mathcal{Z}\left( \mathcal{M%
}\right) \right) :pe\precsim pf\right\} .  \label{eq01}
\end{equation}

Let $a\in LS_{h}\left( \mathcal{M}\right) $ be fixed with spectral measure $%
e^{a}$. For $n\in \mathbb{Z}$, define the central projections $p_{n}$ and $%
q_{n}$ by setting
\begin{equation*}
p_{n}=z\left( e^{a}\left( -\infty ,n\right] ,e^{a}\left( n+1,\infty \right)
\right) ,\ \ \ q_{n}=z\left( e^{a}\left( n+1,\infty \right) ,e^{a}\left(
-\infty ,n\right] \right) .
\end{equation*}%
It may be shown that $p_{n}\vee q_{n}=\mathbf{1}$ for all $n\in \mathbb{Z}$,
that $p_{n}\downarrow _{n}$ and $q_{n}\uparrow _{n}$, and that $%
\bigwedge\nolimits_{n\in \mathbb{Z}}p_{n}=\bigwedge\nolimits_{n\in \mathbb{Z}%
}q_{n}=0$ (see \cite{BdPS}, Lemma 6.12). Using these properties, it follows
easily that there exists a mutually orthogonal sequence $\left\{
z_{n}\right\} _{n\in \mathbb{Z}}$ in $P\left( \mathcal{Z}\left( \mathcal{M}%
\right) \right) $ such that $z_{n}\leq p_{n}q_{n+1}$ for all $n$ and $%
\sum_{n\in \mathbb{Z}}z_{n}=\mathbf{1}$ (see \cite{BdPS}, Corollary 6.14).
The operator $b\in \mathcal{Z}\left( LS_{h}\left( \mathcal{M}\right) \right)
$, appearing in Proposition \ref{Prop02}, is now defined by setting $%
b=\sum_{n\in \mathbb{Z}}\left( n+1\right) z_{n}$. Now we are in a position
to prove the following strengthening of Proposition \ref{Prop02}.

\begin{proposition}
\label{Prop03}If $a\in LS_{h}\left( \mathcal{M}\right) $, then there exists $%
b\in \mathcal{Z}\left( LS_{h}\left( \mathcal{M}\right) \right) $ such that
for each $k\in \mathbb{N}$ there exists a partial isometry $v_{k}\in
\mathcal{M}$ satisfying
\begin{equation}
\left\vert \left[ a,v_{k}\right] \right\vert +2^{-k+1}\mathbf{1}\geq
\left\vert a-b\right\vert .  \label{eq04}
\end{equation}
\end{proposition}

\begin{proof}
Applying the above construction to the operator $2^{k}a$, we obtain for each
$k\in \mathbb{N}$ sequences $\left\{ p_{n}^{k}\right\} _{n\in \mathbb{Z}}$
and $\left\{ q_{n}^{k}\right\} _{n\in \mathbb{Z}}$ of central projections.
We claim that $p_{n}^{k}\leq p_{2n}^{k+1}$. Indeed,
\begin{eqnarray*}
p_{n}^{k}e^{2^{k+1}a}\left( -\infty ,2n\right] &=&p_{n}^{k}e^{a}\left(
-\infty ,n2^{-k}\right] =p_{n}^{k}e^{2^{k}a}\left( -\infty ,n\right] \\
&\precsim &p_{n}^{k}e^{2^{k}a}\left( n+1,\infty \right)
=p_{n}^{k}e^{a}\left( \left( n+1\right) 2^{-k},\infty \right) \\
&\leq &p_{n}^{k}e^{a}\left( \left( 2n+1\right) 2^{-k-1},\infty \right)
=p_{n}^{k}e^{2^{k+1}a}\left( \left( 2n+1\right) ,\infty \right) .
\end{eqnarray*}%
Now, it follows from (\ref{eq01}) and the definition of $p_{2n}^{k+1}$ that $%
p_{n}^{k}\leq p_{2n}^{k+1}$. Via a similar argument it follows also that $%
q_{n}^{k}\leq q_{2n}^{k+1}$.

Using that $p_{n}^{k}\downarrow _{n}$, $q_{n}^{k}\uparrow _{n}$, $%
p_{n}^{k}\leq p_{2n}^{k+1}$ and $q_{n}^{k}\leq q_{2n}^{k+1}$, it follows
easily (using the distributivity of the Boolean algebra $P\left( \mathcal{Z}%
\left( \mathcal{M}\right) \right) $) that
\begin{equation}
p_{n}^{k}q_{n+1}^{k}\leq \left( p_{2n}^{k+1}q_{2n+1}^{k+1}\right) \vee
\left( p_{2n+1}^{k+1}q_{2n+2}^{k+1}\right) .  \label{eq02}
\end{equation}%
This implies that there exist disjoint sequences $\left\{ z_{n}^{k}\right\}
_{n=1}^{\infty }$, $k\in \mathbb{N}$, in $P\left( \mathcal{Z}\left( \mathcal{%
M}\right) \right) $ satisfying $z_{n}^{k}\leq p_{n}^{k}q_{n+1}^{k}$ such
that $z_{n}^{k}=z_{2n}^{k+1}+z_{2n+1}^{k+1}$ for all $n\in \mathbb{Z}$, $%
k\in \mathbb{N}$, and $\sum_{n\in \mathbb{Z}}z_{n}^{k}=\mathbf{1}$ for all $%
k $. Indeed, for $k=1$ this has been observed already above (see \cite{BdPS}%
, Corollary 6.14). If $k\in \mathbb{N}$ is such that $\left\{
z_{n}^{l}\right\} _{n=1}^{\infty }$ has been constructed with the desired
properties for $1\leq l\leq k$, then it follows from (\ref{eq02}) that
\begin{equation*}
z_{n}^{k}\leq \left( p_{2n}^{k+1}q_{2n+1}^{k+1}\right) \vee \left(
p_{2n+1}^{k+1}q_{2n+2}^{k+1}\right) ,\ \ \ n\in \mathbb{Z},
\end{equation*}%
and hence we may write $z_{n}^{k}=z_{2n}^{k+1}+z_{2n+1}^{k+1}$, where $%
z_{2n}^{k+1}$, $z_{2n+1}^{k+1}\in P\left( \mathcal{Z}\left( \mathcal{M}%
\right) \right) $ satisfy $z_{2n}^{k+1}\leq p_{2n}^{k+1}q_{2n+1}^{k+1}$ and $%
z_{2n+1}^{k+1}\leq p_{2n+1}^{k+1}q_{2n+2}^{k+1}$. This proves the claim.

Defining for $k\in \mathbb{N}$ the operator $b_{k}\in \mathcal{Z}\left(
LS_{h}\left( \mathcal{M}\right) \right) $ by
\begin{equation*}
b_{k}=2^{-k}\sum_{n\in \mathbb{Z}}\left( n+1\right) z_{n}^{k},
\end{equation*}%
it follows from Proposition \ref{Prop02} (and the discussion following it)
that there exists a partial isometry $v_{k}\in \mathcal{M}$ such that $%
\left\vert \left[ 2^{k}a,v_{k}\right] \right\vert +\mathbf{1}\geq \left\vert
2^{k}a-2^{k}b_{k}\right\vert $, that is,
\begin{equation}
\left\vert \left[ a,v_{k}\right] \right\vert +2^{-k}\mathbf{1}\geq
\left\vert a-b_{k}\right\vert ,\ \ \ k\in \mathbb{N}.  \label{eq05}
\end{equation}%
Next, we observe that
\begin{eqnarray*}
b_{k}-b_{k+1} &=&\left( \sum\nolimits_{n\in \mathbb{Z}}2^{-k}\left(
n+1\right) z_{2n}^{k+1}+\sum\nolimits_{n\in \mathbb{Z}}2^{-k}\left(
n+1\right) z_{2n+1}^{k+1}\right) \\
&&-\left( \sum\nolimits_{n\in \mathbb{Z}}2^{-k-1}\left( 2n+1\right)
z_{2n}^{k+1}+\sum\nolimits_{n\in \mathbb{Z}}2^{-k-1}\left( 2n+2\right)
z_{2n+1}^{k+1}\right) \\
&=&\sum\nolimits_{n\in \mathbb{Z}}2^{-k-1}z_{2n}^{k+1},
\end{eqnarray*}%
which shows that
\begin{equation}
0\leq b_{k}-b_{k+1}\leq 2^{-k-1}\mathbf{1},\ \ \ k\in \mathbb{N}.
\label{eq03}
\end{equation}%
This estimate implies that the series $c=\sum_{k=1}^{\infty }\left(
b_{k}-b_{k+1}\right) $ is norm convergent in $\mathcal{Z}\left( \mathcal{M}%
\right) $. Defining $b\in \mathcal{Z}\left( LS_{h}\left( \mathcal{M}\right)
\right) $ by setting $b=b_{1}-c$, it follows easily from (\ref{eq03}) that $%
0\leq b-b_{k}\leq 2^{-k}\mathbf{1}$ for all $k\in \mathbb{N}$.

Since $b_{k}$ and $b$ belong to $\mathcal{Z}\left( LS_{h}\left( \mathcal{M}%
\right) \right) $, it is also clear that
\begin{equation*}
\left\vert a-b_{k}\right\vert \geq \left\vert a-b\right\vert -\left(
b-b_{k}\right) \geq \left\vert a-b\right\vert -2^{-k}\mathbf{1}
\end{equation*}%
for all $k$. Therefore, (\ref{eq04}) is now an immediate consequence of (\ref%
{eq05}). The proof is complete.\medskip
\end{proof}

\begin{remark}
To illustrate the relation between operator estimates and norm estimates,
suppose that $a\in LS_{h}\left( \mathcal{M}\right) $ is such that the
derivation $\delta $, defined by $\delta \left( x\right) =\left[ a,x\right] $%
, $x\in LS\left( \mathcal{M}\right) $, satisfies $\delta \left( \mathcal{M}%
\right) \subseteq \mathcal{M}$. If the operator $b\in \mathcal{Z}\left(
LS_{h}\left( \mathcal{M}\right) \right) $ and the partial isometries $%
v_{k}\in \mathcal{M}$ are as in Proposition \ref{Prop03}, then $\left\vert
\delta \left( v_{k}\right) \right\vert +2^{-k+1}\mathbf{1}\geq \left\vert
a-b\right\vert $. This implies that $a-b\in \mathcal{M}$ and $\left\Vert
a-b\right\Vert \leq 2^{-k+1}+\left\Vert \delta \right\Vert $, $k\in \mathbb{N%
}$, and hence, $\left\Vert a-b\right\Vert \leq \left\Vert \delta \right\Vert
$. Evidently, the operator $a-b$ induces the same derivation $\delta $ on $%
LS_{h}\left( \mathcal{M}\right) $ (and on $\mathcal{M}$).
\end{remark}

If $\mathcal{M}$ is a $II_{1}$ factor, then the above result may be further
strengthened. In the proof of the next proposition, we shall make use of the
following facts. The unitary group of $\mathcal{M}$ is denoted by $U\left(
\mathcal{M}\right) $.

\begin{lemma}
\label{Lem02}Let $\mathcal{M}$ be a type $II_{1}$ factor with a fixed normal
faithful finite trace $\tau $.

\begin{enumerate}
\item[(i).] If $e,f\in P\left( \mathcal{M}\right) $ satisfy $\tau \left(
e\right) =\tau \left( f\right) $, then $e\sim f$.

\item[(ii).] If $e,f\in P\left( \mathcal{M}\right) $ are such that $e\sim f$%
, then there exists $u\in U\left( \mathcal{M}\right) $ such that $f=u^{\ast
}eu$. If, in addition, $ef=fe=0$ and $e+f=\mathbf{1}$, then also $%
f=ueu^{\ast }$.
\end{enumerate}
\end{lemma}

\begin{proof}
(i). This follows from \cite{Dix}, Proposition III.2.7.13.

(ii). The first assertion follows from \cite{KR}, Exercise 6.9.11. Assuming,
in addition, that $e+f=\mathbf{1}$, it follows that $ufu^{\ast }=e=\mathbf{1}%
-f=u\left( \mathbf{1}-u^{\ast }fu\right) u^{\ast }$ and so, $f=\mathbf{1}%
-u^{\ast }fu$. This implies that $u^{\ast }fu=\mathbf{1}-f=e$, that is, $%
f=ueu^{\ast }$. \medskip
\end{proof}

\begin{lemma}
\label{Lem03}Let $\mathcal{M}$ be a von Neumann algebra. Suppose that $b\in
LS_{h}\left( \mathcal{M}\right) $ and $p\in P\left( \mathcal{M}\right) $ are
such that $bp=pb$, and define $c=bp-b\left( 1-p\right) $. If $c\geq 0$, then
$c=\left\vert b\right\vert $.
\end{lemma}

\begin{proof}
A simple computation shows that $c^{2}=b^{2}p+b^{2}\left( \mathbf{1}%
-p\right) =b^{2}$ and hence, $c=\left\vert b\right\vert $.
\end{proof}

\begin{proposition}
If $\mathcal{M}$ is a $II_{1}$ factor and $a\in S_{h}\left( \mathcal{M}%
\right) $, then there exist $\lambda _{0}\in \mathbb{R}$ and $u\in U(%
\mathcal{M})$ that $|au-ua|=|u^{\ast }au-a|=u^{\ast }|a-\lambda _{0}\mathbf{1%
}|u+|a-\lambda _{0}\mathbf{1}|$.
\end{proposition}

\begin{proof}
Let $\tau $ be a faithful normal trace on $\mathcal{M}$ satisfying $\tau
\left( \mathbf{1}\right) =1$. It should be observed that $\lim_{\lambda
\rightarrow \infty }\tau \left( e^{a}\left( \lambda ,\infty \right) \right)
=0$ and $\lim_{\lambda \rightarrow -\infty }\tau \left( e^{a}\left( \lambda
,\infty \right) \right) =1$. Define $\lambda _{0}\in \mathbb{R}$ by setting
\begin{equation*}
\lambda _{0}=\inf \left\{ \lambda \in \mathbb{R}:\tau \left( e^{a}\left(
\lambda ,\infty \right) \right) \leq 1/2\right\} .
\end{equation*}%
The normality of the trace implies that $\tau \left( e^{a}\left(
\lambda _{0},\infty \right) \right) \leq 1/2$. Furthermore, if
$\lambda <\lambda _{0} $, then $\tau \left( e^{a}\left( \lambda
,\infty \right) \right) >1/2$ and so, $\tau \left( e^{a}\left(
-\infty ,\lambda \right] \right) <1/2$. This implies that $\tau
\left( e^{a}\left( -\infty ,\lambda _{0}\right) \right) \leq 1/2$.
Since $\mathcal{M}$ does not contain any minimal projections, it
follows that there exist projections $p,q\in P\left(
\mathcal{M}\right)$ such that $e^{a}\left\{ \lambda _{0}\right\}
=p+q,\ pq=0$ and
\begin{equation*}
\tau \left( e^{a}\left( -\infty ,\lambda _{0}\right) +p\right) =\tau \left(
e^{a}\left( \lambda _{0},\infty \right) +q\right) =1/2.
\end{equation*}%
Defining $e,f\in P\left( \mathcal{M}\right) $ by setting $e=e^{a}\left(
-\infty ,\lambda _{0}\right) +p$ and $f=e^{a}\left( \lambda _{0},\infty
\right) +q$, it follows from Lemma \ref{Lem02} that there exists $u\in
U\left( \mathcal{M}\right) $ such that $f=u^{\ast }eu=ueu^{\ast }$ (and so, $%
e=u^{\ast }fu=ufu^{\ast }$).

It should also be observed that the projections $p$ and $q$ commute with the
spectral measure $e^{a}$ and hence, $p$ and $q$ commute with $a$. Since $%
ae^{a}\left\{ \lambda _{0}\right\} =\lambda _{0}e^{a}\left\{ \lambda
_{0}\right\} $, this implies that $ap=\lambda _{0}p$ and $aq=\lambda _{0}q$.
We claim that
\begin{equation*}
\left( a-\lambda _{0}\mathbf{1}\right) \left( f-e\right) =\left\vert
a-\lambda _{0}\mathbf{1}\right\vert .
\end{equation*}%
Indeed,
\begin{eqnarray*}
\left( a-\lambda _{0}\mathbf{1}\right) \left( f-e\right)  &=&\left(
a-\lambda _{0}\mathbf{1}\right) \left[ e^{a}\left( \lambda _{0},\infty
\right) -e^{a}\left( -\infty ,\lambda _{0}\right) +q-p\right]  \\
&=&\left( a-\lambda _{0}\mathbf{1}\right) \left[ e^{a}\left( \lambda
_{0},\infty \right) -e^{a}\left( -\infty ,\lambda _{0}\right) \right]
=\left\vert a-\lambda _{0}\mathbf{1}\right\vert ,
\end{eqnarray*}%
which proves the claim. This implies that
\begin{eqnarray*}
u^{\ast }\left\vert a-\lambda _{0}\mathbf{1}\right\vert u &=&u^{\ast }\left(
a-\lambda _{0}\mathbf{1}\right) uu^{\ast }\left( f-e\right) u \\
&=&\left( u^{\ast }au-\lambda _{0}\mathbf{1}\right) \left( e-f\right) .
\end{eqnarray*}%
Consequently,
\begin{eqnarray*}
\left\vert a-\lambda _{0}\mathbf{1}\right\vert +u^{\ast }\left\vert
a-\lambda _{0}\mathbf{1}\right\vert u &=&\left( a-\lambda _{0}\mathbf{1}%
\right) \left( f-e\right) +\left( u^{\ast }au-\lambda _{0}\mathbf{1}\right)
\left( e-f\right)  \\
&=&a\left( f-e\right) +\left( u^{\ast }au\right) \left( e-f\right)  \\
&=&\left( u^{\ast }au-a\right) e-\left( u^{\ast }au-a\right) f\geq 0.
\end{eqnarray*}%
By an appeal to Lemma \ref{Lem03}, we may conclude that
\begin{equation*}
\left\vert a-\lambda _{0}\mathbf{1}\right\vert +u^{\ast }\left\vert
a-\lambda _{0}\mathbf{1}\right\vert u=\left\vert u^{\ast }au-a\right\vert .
\end{equation*}%
Finally, it is easily verified that $\left\vert u^{\ast }au-a\right\vert
=\left\vert u^{\ast }\left( au-ua\right) \right\vert =\left\vert
au-ua\right\vert $, so the proof is complete.\medskip
\end{proof}

A.F. Ber

Chief software developer

ISV "Solutions"

Tashkent, Uzbekistan

e-mail: ber@ucd.uz

\bigskip

B. de Pagter

Delft Institute of Applied Mathematics

Faculty EEMCS, Delft University of Technology

P.O. Box 5031, 2600 GA Delft, The Netherlands

e-mail: b.depagter@tudelft.nl

\bigskip

F.A. Sukochev

School of Mathematics and Statistics

University of New South Wales

Kensington, NSW 2052, Australia

e-mail: f.sukochev@unsw.edu.au

\end{document}